\documentclass[11pt]{article}
\usepackage{mathrsfs}
\usepackage{epic,eepic,epsf,epsfig}
\usepackage{amsfonts,srcltx,mathrsfs}
\usepackage{amsmath}
\usepackage{amssymb}
\usepackage{amsbsy}
\usepackage{graphicx}
\usepackage{amsfonts}
\usepackage{color}
\usepackage{xcolor}
\usepackage{mathrsfs}
\usepackage{epic,eepic,epsf,epsfig}
\usepackage{graphicx}
\usepackage{amsfonts,srcltx,mathrsfs}
\usepackage{array}
\usepackage{amsthm}
\newtheorem{theorem}{Theorem}
\newtheorem{lem}{Lemma}
\newtheorem{cor}{Corollary}
\newtheorem{defi}{Definition}

\newtheorem{clm}{Claim}
\newtheorem{rmk}{Remark}

\def\f{\noindent}

\begin{document}

\title{Extremal even-cycle-free subgraphs of the $k$-ary $n$-cube}
\author{Shu-Li Zhao\footnote{Corresponding author. Email: zhaoshuli0210@126.com (S.-L. Zhao)}, Yu-Qi Li\\[0.2cm]
{\em\small School of Mathematics and Statistics, Henan Normal University,}\\
{\small\em Xinxiang, China}\\
}

\date{}

\maketitle

\begin{abstract}
A graph is said to be $H$-free if it contains no subgraph isomorphic to $H$. The generalized Tur\'{a}n number $\mathrm{ex}(G, H)$, a natural extension of the classical Tur\'{a}n number, is defined as the maximum number of edges in an $H$-free spanning subgraph of $G$. For any integer $\ell\geq 2$, let $C_{2\ell}$ denote the even cycle of length $2\ell$. So far, extensive studies of the generalized Tur\'{a}n number for the hypercube $Q_{n}$ have been conducted. The $k$-ary $n$-cube is an important family of graphs with many favorable properties, such as regularity and vertex-transitivity. It generalizes the hypercube and reduces to the hypercube when $k=2$. In this paper, we investigate the generalized Tur\'{a}n number $\mathrm{ex}(Q^{k}_{n}, C_{2\ell})$ of $Q^{k}_{n}$ for all even integers $k\geq 4$ and integers $\ell\geq 2$, generalizing the result for $Q_{n}$ [J. Graph Theory 16 (1992) 273-286] from $k=2$ to $k\geq 4$.




\f {\em Keywords:} Tur\'{a}n number; Generalized Tur\'{a}n number; Even-cycle-free subgraph; $k$-ary $n$-cube.
\end{abstract}

\section{Introduction}
Let $H$ be a finite and simple graph. A graph is said to be $H$-free if it contains no subgraph isomorphic to $H$. In 1941, Tur\'{a}n~\cite{Turan-1941} introduced the definition of \textit{Tur\'{a}n number}, which is defined as the maximum number of edges in an $H$-free graph on $n$ vertices and denoted by $\mathrm{ex}(n, H)$.
For the case $H=K_{r+1}$, Tur\'{a}n's theorem gives the exact value of $\mathrm{ex}(n,K_{r+1})$ and characterizes the unique extremal graph, which is the complete $r$-partite graph on $n$ vertices with partitions of nearly equal size.
This result was subsequently extended by Erd\H{o}s and Stone~\cite{ErdosStone-1946} and further refined by Erd\H{o}s and Simonovits~\cite{ErdosSimonovits-1966}, yielding a general asymptotic theorem applicable to arbitrary graphs $H$. For non-bipartite graphs $H$, the theorem of Erd\H{o}s and Simonovits~\cite{ErdosSimonovits-1966} delivers a precise asymptotic formula for $\mathrm{ex}(n,H)$. In contrast, it only gives the trivial bound $\mathrm{ex}(n,H)= o(n^{2})$ for bipartite graphs $H$, demonstrating that the Tur\'{a}n problem is considerably more challenging in the bipartite setting. In the research on this topic, even cycles $C_{2\ell}$ have attracted considerable attention from researchers. A classical result of Bondy and Simonovits~\cite{Bondy-1974} states that for every fixed $\ell\geq 2$, $\mathrm{ex}(n,C_{2\ell})=O(n^{1+\frac{1}{\ell}})$. This bound is known to be tight for $\ell= 2,3,5$.

Later, Erd\H{o}s~\cite{Erdos-1984} extended the {\it Tur\'{a}n number} to the {\it generalized Tur\'{a}n number}, where the host graph is fixed to a given graph $G$, instead of ranging over the family of all $n$-vertex graphs. The {\it generalized Tur\'{a}n number}, denoted by $\mathrm{ex}(G, H)$, is defined as the maximum number of edges in an $H$-free spanning subgraph of $G$. When $H$ is bipartite, the problem is particularly challenging, and even the asymptotic order of $\mathrm{ex}(G, H)$ is often undetermined~\cite{Erdos-1966}. The hypercube $Q_{n}$ serves as a well-studied host graph, and the generalized Tur\'{a}n number $\mathrm{ex}(Q_{n},C_{2\ell})$ for the $n$-dimensional hypercube has attracted considerable attention.
In 1984, Erd\H{o}s~\cite{Erdos-1984} conjectured that $\mathrm{ex}(Q_{n},C_{4})\leq \big(\frac{1}{2}+o(1)\big)e(Q_{n})$, and further raised the question of whether $o(e(Q_{n}))$ edges can force the presence of a $C_{2\ell}$ for all integers $\ell\geq 3$. For the $C_4$-free problem on hypercubes, the first upper bound on $\mathrm{ex}(Q_{n},C_{4})$ was established by Chung~\cite{Chung-1992}. Shortly afterward, Brass, Harborth, and Nienborg~\cite{Brass-1995} derived a corresponding lower bound, proving that $\mathrm{ex}(Q_{n},C_{4})\geq \big(\frac{1}{2}+\frac{c}{\sqrt{n}}\big)e(Q_{n})$ for some constant $c>0$. Notably, a substantial gap persisted between the available upper and lower bounds. To narrow this gap, Thomason and Wagner~\cite{Thomason-2009} improved the upper bound in 2009, which was further refined by Balogh, Hu, Lidick\'{y} and Liu~\cite{Balogh-2014} in $2014$ to the value of $(0.6068+o(1))e(Q_{n})$.
For $C_6$-free problem on hypercubes, Chung~\cite{Chung-1992} established the upper bound $\mathrm{ex}(Q_{n},C_{6})\leq \big(\sqrt{2}-1+o(1)\big)e(Q_{n})$. Shortly afterward, Conder~\cite{Conder-1993} derived the corresponding lower bound $\mathrm{ex}(Q_{n},C_{6})\geq \frac{1}{3}e(Q_{n})$. The upper bound was further improved by Balogh, Hu, Lidick\'{y} and Liu~\cite{Balogh-2014}, who proved that $\mathrm{ex}(Q_{n},C_{6})\leq (0.3755+o(1))e(Q_{n})$. Regarding longer even cycles, Chung~\cite{Chung-1992} showed that $\mathrm{ex}(Q_{n},C_{4\ell})=o(e(Q_{n}))$ for all integers $\ell\geq 2$. F\"{u}redi and \"{O}zkahya~\cite{Furedi-2009} subsequently established the result for $C_{4\ell+2}$ with $\ell\geq 3$. A unified approach to even-cycle forbidden problems was subsequently developed by Conlon~\cite{Conlon-2010}, who verified that $\mathrm{ex}(Q_{n},C_{2\ell})=o(e(Q_{n}))$ for all $\ell\geq 4$ with $\ell\neq5$. Beyond hypercubes, generalized Tur\'{a}n numbers have also been studied for Cayley graphs, owing to their desirable structural properties such as regularity and vertex-transitivity. For example, Cao et al.~\cite{Cao-2021} investigated the extremal problem on the complete transposition graph $CT_{n}$. Additionally, Cao et al.~\cite{Cao-2021(1)} considered the analogous extremal problem on doubled Johnson graphs.

The $k$-ary $n$-cube $Q^{k}_{n}$, introduced by Dally~\cite{Dally-1990}, is an important family of graphs. It generalizes the hypercube and reduces to the hypercube when $k=2$~\cite{Mao-2002}. For even $k$, the $k$-ary $n$-cube $Q^{k}_{n}$ is bipartite~\cite{Wang-2014} and has many favorable properties, such as regularity, vertex-transitivity and a recursive structure~\cite{Becker-1988,Ghozati-1999,Mao-2002,Xiang-2011}. While the generalized Tur\'{a}n number $\mathrm{ex}(Q_{n},C_{2\ell})$ has been extensively studied on the hypercube, far less is known about the $k$-ary $n$-cube $Q^{k}_{n}$. In this paper, we mainly study the generalized Tur\'{a}n number $\mathrm{ex}(Q^{k}_{n},C_{2\ell})$ for $Q^{k}_{n}$, where $k$ is even with $k\geq 4$ and $\ell\geq 2$, and we consider the asymptotic regime as $n\rightarrow\infty$.

The remainder of this paper is organized as follows. Section~\ref{sec:Prel} introduces notation, definitions, and preliminary properties of $Q^{k}_{n}$.
In Section~\ref{sec:$2l$-Cycle-free and $l$ is even}, we establish an upper bound on $\mathrm{ex}(Q^{k}_{n},C_{2\ell})$ for even $\ell$ with $\ell\geq 4$.
In Subsection~\ref{sec:6-Cycle-free}, a general upper bound on $\mathrm{ex}(Q^{k}_{n}, C_6)$ is proved and then refined by a more careful argument. In Subsection~\ref{sec:4-Cycle-free}, we establish an upper bound on $\mathrm{ex}(Q^{k}_{n},C_{4})$ for all even $k$ with $k\geq 4$, and improve this bound in the special case $k=4$.
In Section~\ref{sec:remarks}, the paper is concluded.

\section{Preliminary}
\label{sec:Prel}
All graphs considered in this paper are finite and simple. For a graph $G$, let $V(G)$ and $E(G)$ denote its vertex set and edge set, respectively, and let $v(G)=|V(G)|$ and $e(G)=|E(G)|$. For a vertex $v\in V(G)$, let $N_{G}(v)$ denote the set of neighbors of $v$ in $G$, and let $d_{G}(v)=|N_{G}(v)|$ be the degree of $v$. For a positive integer $\ell$, let $C_{\ell}$ denote a cycle of length $\ell$. For subgraphs $H_{1}$ and $H_{2}$ of $G$, define $H_{1}\cap H_{2}$ as the graph with vertex set $V(H_{1})\cap V(H_{2})$ and edge set $E(H_{1})\cap E(H_{2})$. The subgraph induced by $V^{\prime}$ in $G$, denoted by $G[V^{\prime}]$, is a graph whose vertex set is $V^{\prime}$ and the edge set is the set of all the edges of $G$ with both ends in $V^{\prime}$. For any positive integer $n$, let $[n]$ denote the set $\{1,2,\ldots,n\}$ and $\langle n\rangle$ denote the set $\{0,1,2,\ldots,n\}$.

\begin{defi}[\cite{Wang-2015}]\label{defi1}
The $k$-ary $n$-cube, denoted by $Q^{k}_{n}$, with $k\geq 2$ and $n\geq 1$, is a simple graph on $k^{n}$ vertices, each of which is a $k$-ary $n$-dimensional string $u=u_{1}u_{2}\ldots u_{n}$ with $u_{i}\in\langle k-1\rangle$ for each $i\in[n]$. Two vertices $u=u_{1}u_{2}\ldots u_{n}$ and $v=v_{1}v_{2}\ldots v_{n}$ are adjacent if and only if there exists an integer $j\in[n]$ such that $u_{j}=v_{j}\pm1$ (mod $k$) and $u_{i}=v_{i}$ for every $i\in[n]\setminus\{j\}$.
\end{defi}

By Definition~\ref{defi1}, one can obtain that $Q^{k}_{n}$ is $n$-regular for $k=2$ and $2n$-regular for $k\geq 3$. In particular, $Q^{k}_{n}$ is exactly the $n$-dimensional hypercube $Q_{n}$ for $k=2$. Moreover, $Q^{k}_{n}$ is a bipartite $2n$-regular graph~\cite{Wang-2014} with $k^{n}$ vertices and $nk^{n}$ edges for even integer $k$ with $k\geq 4$.

The following lemma is about the girth of $Q^{k}_{n}$.

\begin{lem}[\cite{Li-2011}]\label{lem1}
The girth of $Q^{k}_{n}$ is $4$ for all integers $k\geq 4$ and $n\geq 2$.
\end{lem}

The following two lemmas give the total number of $4$-cycles in $Q^{k}_{n}$, together with the number of $4$-cycles containing a prescribed edge, separately for the case $k\geq 6$ and $k=4$.

\begin{lem}\label{lem2} Let $Q^{k}_{n}$ be the $k$-ary $n$-cube, and let $n_{k}(C_{4})$ denote the number of $4$-cycles in $Q^{k}_{n}$, where $k$ is even and $k\geq 6$. Then the following statements hold.
\begin{enumerate}
\item[{\rm (1)}] Every edge $e\in E(Q^{k}_{n})$ is contained in exactly $2(n-1)$ $4$-cycles of $Q^{k}_{n}$.
\item[{\rm (2)}] $n_{k}(C_{4})=\binom{n}{2}k^{n}$.
\end{enumerate}
\end{lem}

\noindent{\bf Proof.}
(1) Without loss of generality, let $e=uv\in E(Q^{k}_{n})$, where $u=u_{1}\ldots u_{i}\ldots u_{n}$ and $v=u_{1}\ldots (u_{i}+1)\ldots u_{n}$. Take $w\in N_{Q^{k}_{n}}(u)\setminus\{v\}$ such that $u$ and $w$ differ in the $j$-th coordinate. Then $(v,u,w)$ is a $2$-path of $Q^{k}_{n}$.

\begin{clm}\label{clm1} If the $2$-path $(v,u,w)$ is contained in a $4$-cycle, then $j\neq i$. \end{clm}

\noindent{\bf Proof of Claim~\ref{clm1}.} Suppose to the contrary that $j=i$. Since $w\neq v$, by Definition~\ref{defi1}, we have \[ w=u_{1}\ldots (u_{i}-1)\ldots u_{n}. \] Assume that $(v,u,w)$ is contained in a $4$-cycle $(v,u,w,x,v)$, and suppose that $w$ and $x$ differ in the $s$-th coordinate. If $s=i$, by Definition~\ref{defi1}, \[ x=u_{1}\ldots (u_{i}-2)\ldots u_{n}. \] Hence $x$ and $v$ differ only in the $i$-th coordinate, where they differ by exactly $3$. Thus, $xv\notin E(Q^{k}_{n})$, a contradiction. If $s\neq i$, then $x$ differs from $v$ in both the $i$-th and the $s$-th coordinates. By Definition~\ref{defi1}, $xv\notin E(Q^{k}_{n})$, a contradiction.  Therefore, $j\neq i$. \hfill$\blacksquare$ \vskip0.2cm

Now suppose that $j\neq i$. Without loss of generality, assume $1\leq i<j\leq n$ and denote $w=u_{1}\ldots u_{i}\ldots (u_{j}+\varepsilon)\ldots u_{n}$ for some $\varepsilon\in\{+1,-1\}$. Take $x\in N_{Q^{k}_{n}}(w)\setminus\{u\}$ such that $w$ and $x$ differ in the $s$-th coordinate. Then $(v,u,w, x)$ is a $3$-path in $Q^{k}_{n}$.

\begin{clm}\label{clm2} If the $3$-path $(v,u,w,x)$ is contained in the $4$-cycle $(v,u,w,x,v)$, then $s=i$.
\end{clm}

\noindent{\bf Proof of Claim~\ref{clm2}.} Suppose to the contrary that $s\neq i$, by Definition~\ref{defi1}, $s=j$ or $s\in [n]\setminus\{i,j\}$. If $s=j$, by Definition~\ref{defi1}, \[ x=u_{1}\ldots u_{i}\ldots (u_{j}+2\varepsilon)\ldots u_{n}. \] Then $x$ and $v$ differ in both the $i$-th and the $j$-th coordinates. Hence $xv\notin E(Q^{k}_{n})$, a contradiction. If $s\in [n]\setminus\{i,j\}$, then $x$ and $v$ differ in the $i$-th, the $j$-th and the $s$-th coordinates. Hence $xv\notin E(Q^{k}_{n})$, a contradiction. Thus, $s=i$. \hfill$\blacksquare$ \vskip0.2cm

\begin{clm}\label{clm3} If $j\neq i$ and $\varepsilon\in\{+1,-1\}$, then the $2$-path $(v,u,w)$ is contained in a unique $4$-cycle $(v,u,w,x,v)$, where \[ w=u_{1}\ldots u_{i}\ldots (u_{j}+\varepsilon)\ldots u_{n}, \qquad x=u_{1}\ldots (u_{i}+1)\ldots (u_{j}+\varepsilon)\ldots u_{n}. \] \end{clm}

\noindent{\bf Proof of Claim~\ref{clm3}.} By Claim~\ref{clm2}, $w$ and $x$ must differ in the $i$-th coordinate. Hence there are two possible candidates for $x$, namely \[ x'=u_{1}\ldots (u_{i}+1)\ldots (u_{j}+\varepsilon)\ldots u_{n} \] or \[ x''=u_{1}\ldots (u_{i}-1)\ldots (u_{j}+\varepsilon)\ldots u_{n}. \]
Clearly, $x'$ and $v$ differ only in the $j$-th coordinate, where they differ by $1$. Thus, $x'v\in E(Q^{k}_{n})$. As $x''$ and $v$ differ in both the $i$-th and the $j$-th coordinates, $x''v\notin E(Q^{k}_{n})$. Thus, $x=x'$ and $(v,u,w,x',v)$ is the unique $4$-cycle containing the $2$-path $(v,u,w)$. \hfill$\blacksquare$ \vskip0.2cm

For any given edge $e=uv$, there are $(n-1)$ choices for $j$, namely every $j\in [n]\setminus\{i\}$. For each such $j$, there are two choices $\varepsilon\in\{+1,-1\}$. By Claim~\ref{clm3}, each choice gives rise to a unique $4$-cycle containing $e$. Moreover, these $4$-cycles are pairwise distinct. Consequently, every edge of $Q^{k}_{n}$ is contained in exactly $2(n-1)$ $4$-cycles. \medskip

(2) By part (1), we obtain that every edge is contained in exactly $2(n-1)$ $4$-cycles. Double-counting the pairs $(e, C)$ for which $e$ is an edge of the $4$-cycle $C$ yields
\[ 4n_{k}(C_{4}) =2(n-1)e(Q^{k}_{n}). \]
Hence
\[ n_{k}(C_{4})
=\binom{n}{2}k^{n}. \]
\hfill $\blacksquare$
\vskip0.2cm

\begin{lem}\label{lem3} Let $Q^{4}_{n}$ be the $4$-ary $n$-cube, and let $n_{4}(C_{4})$ denote the number of $4$-cycles in $Q^{4}_{n}$. Then the following statements hold.
\begin{enumerate}
\item[{\rm (1)}] Every edge $e\in E(Q^{4}_{n})$ is contained in exactly $(2n-1)$ $4$-cycles of $Q^{4}_{n}$.
\item[{\rm (2)}] $n_{4}(C_{4})=(2n^{2}-n)4^{n-1}$.
\end{enumerate}
\end{lem}
\noindent{\bf Proof.}
(1) Without loss of generality, let $e=uv\in E(Q^{4}_{n})$, where $u=u_{1}\ldots u_{i}\ldots u_{n}$ and $v=u_{1}\ldots (u_{i}+1)\ldots u_{n}$. Take $w\in N_{Q^{4}_{n}}(u)\setminus\{v\}$ and assume that $u$ and $w$ differ in the $j$-th coordinate. Then $(v,u,w)$ is a $2$-path of $Q^{4}_{n}$.

If $j\neq i$, by Claim~\ref{clm3} of Lemma~\ref{lem2}, the $2$-path $(v,u,w)$ is contained in a unique $4$-cycle. Now suppose that $j=i$. Since $w\neq v$, by Definition~\ref{defi1}, we have
\[ w=u_{1}\ldots (u_{i}-1)\ldots u_{n}. \]
Take $x\in N_{Q^{4}_{n}}(w)\setminus\{u\}$ such that $w$ and $x$ differ in the $s$-th coordinate. Then $(v,u,w,x)$ is a $3$-path in $Q^{4}_{n}$.

\begin{clm}\label{clm4} If the $3$-path $(v,u,w,x)$ is contained in the $4$-cycle $(v,u,w,x,v)$, then $s=i$. \end{clm}

\noindent{\bf Proof of Claim~\ref{clm4}.} Suppose to the contrary that $s\neq i$. Then $x$ differs from $v$ in both the $i$-th and the $s$-th coordinates, and hence $xv\notin E(Q^{4}_{n})$, a contradiction. Thus, $s=i$.
\hfill$\blacksquare$
\vskip0.2cm

\begin{clm}\label{clm5} If $j=i$, then the $2$-path $(v,u,w)$ is contained in a unique $4$-cycle $(v,u,w,x,v)$, where \[ x=u_{1}\ldots (u_{i}-2)\ldots u_{n}. \]
\end{clm}

\noindent{\bf Proof of Claim~\ref{clm5}.} By Claim~\ref{clm4}, $w$ and $x$ must differ in the $i$-th coordinate. Since $x\neq w$, by Definition~\ref{defi1}, \[ x=u_{1}\ldots (u_{i}-2)\ldots u_{n}. \]
Then $x$ and $v$ differ only in the $i$-th coordinate, where they differ by $1$. Thus $(v,u,w,x,v)$ is the unique $4$-cycle containing the $2$-path $(v,u,w)$.
\hfill$\blacksquare$
\vskip0.2cm

Consequently, by Lemma~\ref{lem2}(1) and Claim~\ref{clm5}, we obtain that every edge of $Q^{4}_{n}$ is contained in exactly $(2n-1)$ $4$-cycles, where $2n-1=2(n-1)+1$. \medskip

(2) By part (1), we obtain that every edge is contained in exactly $(2n-1)$ $4$-cycles. Double-counting the pairs $(e,C)$ for which $e$ is an edge of the $4$-cycle $C$ gives
\[ 4n_{4}(C_{4}) =(2n-1)e(Q^{4}_{n}). \]
Hence
\[ n_{4}(C_{4})
=(2n^{2}-n)4^{n-1}. \]
\hfill$\blacksquare$
\vskip0.2cm

The following corollary is an immediate consequence of Lemma~\ref{lem2} and Lemma~\ref{lem3}. It concerns the existence and uniqueness of a 4-cycle containing a given 2-path in $Q^{k}_{n}$.

\begin{cor}\label{cor1}
Let $(v,u,w)$ be a $2$-path in $Q^{k}_{n}$. It is contained in at most one $4$-cycle for $k\geq6$ and exactly one $4$-cycle for $k=4$.
\end{cor}
\f {\bf Proof.}
Let $(v,u,w)$ be a $2$-path in $Q^{k}_{n}$. Suppose $u$ and $v$ differ in the $i$-th dimension, and $u$ and $w$ differ in the $j$-th dimension, where $i,j\in[n]$. By Lemma~\ref{lem2} and Lemma~\ref{lem3}, if $j\neq i$, then the $2$-path $(v,u,w)$ is contained in a unique $4$-cycle for all even $k\ge 4$. By Claim~\ref{clm1} and Claim~\ref{clm5}, if $j=i$, then the $2$-path $(v,u,w)$ is not contained in any $4$-cycle for $k\ge 6$, and is contained in a unique $4$-cycle for $k=4$.

\hfill $\blacksquare$
\vskip0.2cm

\section{$2\ell$-Cycle-free subgraph of $Q^{k}_{n}$ for even $\ell\geq 4$}
\label{sec:$2l$-Cycle-free and $l$ is even}

In this section, we focus on the case of even integers $\ell\geq 4$. Let $\tilde{G}$ denote a $2\ell$-cycle-free spanning subgraph of $Q^{k}_{n}$ with maximum number of edges. Then, $\mathrm{ex}(Q^{k}_{n},C_{2\ell})=e(\tilde{G})$. To estimate $e(\tilde{G})$, we introduce the following auxiliary graph $G_{v}$.

\begin{defi}\label{defi2}
For each $v\in V(Q^{k}_{n})$, the auxiliary graph $G_{v}$ is the graph with vertex set $V(G_{v})=N_{Q^{k}_{n}}(v)$. Two vertices $u_{1},u_{2}\in V(G_{v})$ are adjacent in $G_{v}$ if and only if they have a common neighbor $w$ in $\tilde{G}$, where $w\neq v$.
\end{defi}

By Definition~\ref{defi2}, we have
\begin{equation}
\sum_{v\in V(Q^{k}_{n})}v(G_{v})=2nk^{n}.
\end{equation}
By Corollary~\ref{cor1}, we obtain
\begin{equation}
\sum_{v\in V(Q^{4}_{n})}e(G_{v})=\sum_{w\in V(Q^{4}_{n})}\binom{d_{\tilde{G}}(w)}{2} \quad \text{for } k=4.
\end{equation}
By Corollary~\ref{cor1}, for any $w\in V(Q^{k}_{n})$ with $k\geq 6$, there exist $n$ $2$-paths centered at $w$ that are not contained in any $4$-cycle. Thus, we obtain
\begin{equation}
\sum_{v\in V(Q^{k}_{n})}e(G_{v})\geq \sum_{w\in V(Q^{k}_{n})}\binom{d_{\tilde{G}}(w)}{2}-nk^{n} \quad \text{for even } k\geq 6.
\end{equation}
Combining (2) with (3), we conclude that
\begin{equation}
\sum_{v\in V(Q^{k}_{n})}e(G_{v})\geq \sum_{w\in V(Q^{k}_{n})}\binom{d_{\tilde{G}}(w)}{2}-nk^{n} \quad \text{for even } k\geq 4.
\end{equation}

The following lemma, due to Bondy and Simonovits~\cite{Bondy-1974}, is useful for our proof.
\begin{lem}[\cite{Bondy-1974}]\label{lem4}
Let $G$ be a graph on $n$ vertices. If $e(G)>c_{k}n^{1+\frac{1}{k}}$, where $c_{k}$ is a positive constant depending only on $k$, then $G$ contains an even cycle $C_{2s}$ for every integer $k\leq s\leq kn^{\frac{1}{k}}$.
\end{lem}

The following theorem gives an upper bound on the generalized Tur\'{a}n number $\mathrm{ex}(Q^{k}_{n},C_{2\ell})$, where $\ell\geq 4$ and $\ell$ is even.

\begin{theorem}\label{theorem1}
If $k$ and $\ell$ are even integers with $k\geq 4$ and $\ell\geq 4$, then $\mathrm{ex}(Q^{k}_{n},C_{2\ell})=O(n^{\frac{1}{\ell}-\frac{1}{2}})e(Q^{k}_{n}).$
\end{theorem}
\f {\bf Proof.}
On the one hand, since $\tilde{G}$ is $C_{2\ell}$-free and $\ell$ is even, $G_{v}$ contains no $\ell$-cycles. Applying Lemma~\ref{lem4} with $s=\frac{\ell}{2}$, we obtain \begin{equation}
e(G_{v})\leq c_{\ell}(v(G_{v}))^{1+\frac{2}{\ell}},
\end{equation}
where $c_{\ell}$ is a positive constant depending only on $\ell$.\\
Combining (5) with (1), we have
\begin{equation}
\sum_{v\in V(Q^{k}_{n})}e(G_{v})\leq c_{\ell}(2n)^{1+\frac{2}{\ell}}k^{n}.
\end{equation}

On the other hand, let $\alpha=\frac{e(\tilde{G})}{e(Q^{k}_{n})}=\frac{e(\tilde{G})}{nk^{n}}$. By the Cauchy-Schwarz inequality and the identity $\sum_{w\in V(Q^{k}_{n})}d_{\tilde{G}}(w)=2e(\tilde{G})$, it follows that
\begin{align}
\displaystyle\sum\limits_{w\in V(Q^{k}_{n})}\binom{d_{\tilde{G}}(w)}{2}
& =\frac{1}{2}\sum\limits_{w\in V(Q^{k}_{n})}d_{\tilde{G}}^{2}(w)-\dfrac{1}{2}\sum\limits_{w\in V(Q^{k}_{n})}d_{\tilde{G}}(w) \notag \\[2pt]
& \geq \dfrac{1}{2k^{n}}\left(\displaystyle\sum\limits_{w\in V(Q^{k}_{n})}d_{\tilde{G}}(w)\right)^{2}-\dfrac{1}{2}\displaystyle\sum\limits_{w\in V(Q^{k}_{n})}d_{\tilde{G}}(w) \notag \\[4pt]
& =2\alpha^{2}n^{2}k^{n}-\alpha nk^{n}.
\end{align}
Combining (4) with (7) yields
\begin{align}
\sum_{v\in V(Q_n^k)} e(G_v)
\geq 2\alpha^2 n^2 k^n - (1+\alpha)n k^n. \label{eq:final}
\end{align}

By (6) and (8), we obtain
$$2\alpha^{2}n^{2}k^{n}-(1+\alpha)nk^{n}\leq c_{\ell}(2n)^{1+\frac{2}{\ell}}k^{n}.$$
Dividing both sides by $nk^{n}$, we get
\begin{align*}
2\alpha^2 n - (1+\alpha) &\le c_\ell 2^{1+\frac{2}{\ell}} n^{\frac{2}{\ell}} \\[2pt]
\alpha^2 &\le c_\ell 2^{\frac{2}{\ell}} n^{\frac{2}{\ell} - 1} + \frac{\alpha+1}{2} n^{-1} \\
\alpha^2 &\le \left(c_\ell 2^{\frac{2}{\ell}} + \frac{\alpha+1}{2}\right) n^{\frac{2}{\ell} - 1}.
\end{align*}
Let $c=\sqrt{c_{\ell}2^{\frac{2}{\ell}}+\frac{\alpha+1}{2}}$, which is a positive constant depending only on $\ell$. Then,

$$\alpha\leq cn^{\frac{1}{\ell}-\frac{1}{2}}.$$
As $\alpha=\frac{e(\tilde{G})}{e(Q^{k}_{n})}$, we obtain that
$$\mathrm{ex}(Q^{k}_{n},C_{2\ell})=e(\tilde{G})=\alpha e(Q^{k}_{n})=O(n^{\frac{1}{\ell}-\frac{1}{2}})e(Q^{k}_{n}).$$
\hfill $\blacksquare$
\vskip0.2cm

Recall that $\alpha=\frac{e(\tilde{G})}{e(Q^{k}_{n})}$. By (7), we have
\begin{align*}
\sum_{w\in V(Q_n^k)} \binom{d_{\tilde{G}}(w)}{2} - nk^n
\geq nk^n (2n\alpha^2 - \alpha - 1).
\end{align*}
It follows that the left-hand side is nonnegative whenever $\alpha\geq \frac{1+\sqrt{1+8n}}{4n}$. Thus, we have the following remark.

\begin{rmk}
The inequality $\sum_{w\in V(Q^{k}_{n})}\binom{d_{\tilde{G}}(w)}{2}-nk^{n}\geq 0$ holds whenever $e(\tilde{G})=\Omega(n^{-\frac{1}{2}})e(Q^{k}_{n})$.
\end{rmk}
\vskip0.2cm

\section{$2\ell$-Cycle-free subgraph of $Q^{k}_{n}$ for $\ell$=2,3}
\label{sec:6,4-Cycle-free}
Let $\tilde{G}$ denote a $2\ell$-cycle-free spanning subgraph of $Q^{k}_{n}$ with maximum number of edges, where $\ell=2,3$. The following lemma plays a key role in the proof of the main theorem. Recall that $C\cap \tilde{G}$ denotes the graph with vertex set $V(C)\cap V(\tilde{G})$ and edge set $E(C)\cap E(\tilde{G})$.

\begin{figure}[h!]
\begin{center}
\includegraphics[height=8cm,width=12cm]{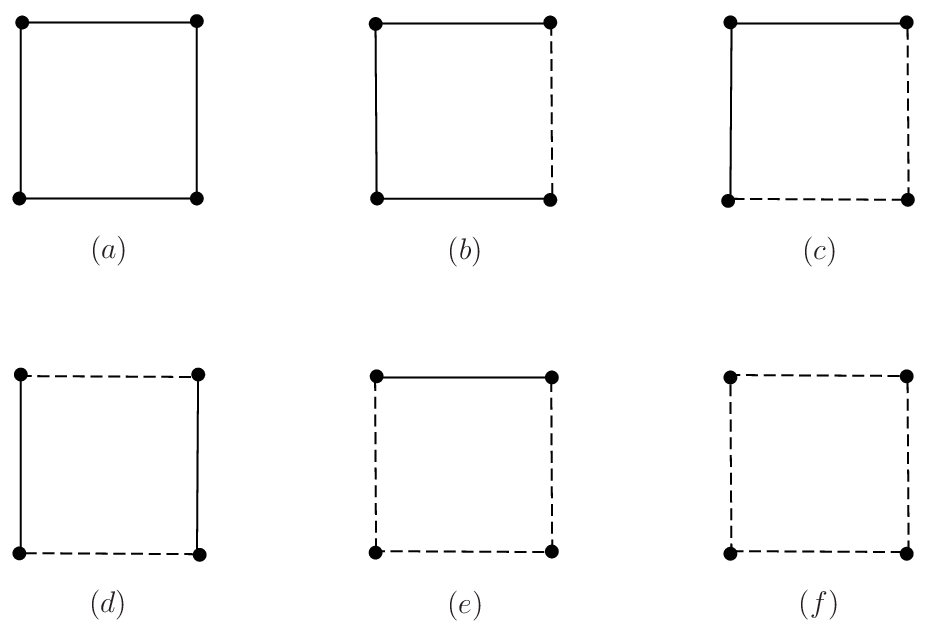}
\end{center}
\begin{center}
\vspace{-0.6cm}
\caption{ Possibilities for $C\cap\tilde{G}$ when $C$ is any 4-cycle in $Q^{k}_{n}$.}\label{F1}
\vspace{-0.6cm}
\end{center}
\end{figure}

\begin{lem}\label{lem5}
Let $\alpha=\frac{e(\tilde{G})}{e(Q^{k}_{n})}$ and let $\theta_{4},\theta_{3},\theta^{1}_{2}, \theta^{2}_{2},\theta_{1},\theta_{0}$ denote the proportions of $4$-cycles $C$ in $Q^{k}_{n}$ with $C\cap \tilde{G}$ isomorphic to the graphs $(a)-(f)$ in Figure~\ref{F1}, respectively. Then
\begin{equation}
\theta_{4}+\theta_{3}+\theta^{1}_{2}+\theta^{2}_{2}+\theta_{1}+\theta_{0}=1,
\end{equation}
and
\begin{equation}
4\theta_{4}+3\theta_{3}+2\theta^{1}_{2}+2\theta^{2}_{2}+\theta_{1}=4\alpha.
\end{equation}
\end{lem}
\f {\bf Proof.}
For any 4-cycle $C$ in $Q^{k}_{n}$, $C\cap \tilde{G}$ is isomorphic to one of the six graphs in Figure~\ref{F1}. Hence
$$\theta_{4}+\theta_{3}+\theta^{1}_{2}+\theta^{2}_{2}+\theta_{1}+\theta_{0}=1.$$

Let $n_{k}(C_{4})$ denote the number of $4$-cycles in $Q^{k}_{n}$ for even $k\geq 4$. For $k\geq 6$, double-counting the pairs $(e,C)$ such that $C$ is a $4$-cycle in $Q^{k}_{n}$ and $e \in E(C \cap \tilde{G})$ yields
$$(4\theta_{4}+3\theta_{3}+2\theta^{1}_{2}+2\theta^{2}_{2}+\theta_{1})n_{k}(C_{4})=2(n-1)e(\tilde{G}).$$
By Lemma~\ref{lem2}(2) and $e(\tilde{G})=\alpha e(Q^{k}_{n})=\alpha nk^{n}$, we obtain
$$(4\theta_{4}+3\theta_{3}+2\theta^{1}_{2}+2\theta^{2}_{2}+\theta_{1})\binom{n}{2}k^{n}= 2\alpha k^{n}n(n-1).$$
Therefore,
\begin{equation}
4\theta_{4}+3\theta_{3}+2\theta^{1}_{2}+2\theta^{2}_{2}+\theta_{1}= 4\alpha.
\end{equation}

For $k=4$, by double-counting the pairs $(e, C)$ for which $C$ is a $4$-cycle in $Q^{4}_{n}$ and $e \in E(C \cap \tilde{G})$, we obtain
$$(4\theta_{4}+3\theta_{3}+2\theta^{1}_{2}+2\theta^{2}_{2}+\theta_{1})n_{4}(C_{4})
=(2n-1)e(\tilde{G}).$$
As $e(\tilde{G}) = \alpha n4^n$, by Lemma~\ref{lem3}(2), we obtain
\begin{equation}
4\theta_{4}+3\theta_{3}+2\theta^{1}_{2}+2\theta^{2}_{2}+\theta_{1}= 4\alpha.
\end{equation}


Equations (11) and (12) together imply that
$$4\theta_{4}+3\theta_{3}+2\theta^{1}_{2}+2\theta^{2}_{2}+\theta_{1}=4\alpha$$
for even $k\geq 4$.
\hfill $\blacksquare$
\vskip0.2cm

\subsection{6-Cycle-free subgraph of $Q^{k}_{n}$}
\label{sec:6-Cycle-free}
In this section, we focus on the case $\ell=3$. Let $\tilde{G}$ denote a 6-cycle-free spanning subgraph of $Q^{k}_{n}$ with maximum number of edges. Recall that $C\cap \tilde{G}$ is isomorphic to one of the graphs $(a)-(f)$ in Figure~\ref{F1} for any 4-cycle $C$ in $Q^{k}_{n}$.

\begin{lem}\label{lem6}
Let $\theta_{3}$ and $\theta_{4}$ denote the proportions of $4$-cycles $C$ in $Q^{k}_{n}$ such that $C\cap \tilde{G}$ contains exactly three edges and four edges, respectively. Then, $\theta_{4}=o(1)$ and $\theta_{3}=o(1)$.
\end{lem}
\f {\bf Proof.}
For any $e\in E(Q^{k}_{n})$, we claim that there is at most one $4$-cycle $C$ in $Q^{k}_{n}$ containing $e$ such that $C-e$ contains three edges of $\tilde{G}$. Suppose to the contrary that there exist two distinct 4-cycles $C^{\prime}$ and $C^{\prime\prime}$ in $Q^{k}_{n}$ containing $e$. By Corollary~\ref{cor1}, $C^{\prime}$ and $C^{\prime\prime}$ share exactly one common edge $e$. Thus, $(C^{\prime}-e)\cup(C^{\prime\prime}-e)$ forms a 6-cycle in $\tilde{G}$, contradicting the assumption that $\tilde{G}$ is 6-cycle-free.

Let $n_{k}(C_{4})$ denote the number of $4$-cycles in $Q^{k}_{n}$ for even $k\geq 4$. If $k\geq 6$, by the claim above, we have
$$(4\theta_{4}+\theta_{3})n_{k}(C_{4})\leq e(Q^{k}_{n}).$$
As $e(Q^{k}_{n})=nk^{n}$, by Lemma~\ref{lem2}(2), we obtain
$$(4\theta_{4}+\theta_{3})\binom{n}{2}k^{n} \leq nk^{n}.$$
Dividing both sides by $\binom{n}{2}k^{n}$ yields
$$4\theta_{4}+\theta_{3} \leq \dfrac{2}{n-1}.$$
Since $n\rightarrow\infty$, we obtain $\theta_{4}=o(1)$ and $\theta_{3}=o(1).$

If $k=4$, the claim above gives
$$(4\theta_{4}+\theta_{3})n_{4}(C_{4})\leq e(Q^{4}_{n}).$$
As $e(Q^{4}_{n})=n4^{n}$, by Lemma~\ref{lem3}(2), we obtain
$$4\theta_{4}+\theta_{3} \leq \dfrac{4}{2n-1}.$$
Since $n\rightarrow\infty$, we obtain $\theta_{4}=o(1)$ and $\theta_{3}=o(1)$.

Combining the two cases $k=4$ and $k\geq 6$, we conclude that
$$\theta_{4}=o(1) \quad \text{and} \quad \theta_{3}=o(1)$$
for all even $k\geq 4$.
\hfill $\blacksquare$
\vskip0.2cm

We now derive an upper bound for the generalized Tur\'{a}n number $\mathrm{ex}(Q^{k}_{n},C_{6})$.

\begin{theorem}\label{theorem2}
If $k$ is an even integer with $k\geq 4$, then $\mathrm{ex}(Q^{k}_{n},C_{6})\leq (\frac{1}{2}+o(1))e(Q^{k}_{n}).$
\end{theorem}
\f {\bf Proof.}
Let $\alpha=\frac{e(\tilde{G})}{e(Q^{k}_{n})}$. By Lemma~\ref{lem5} and Lemma~\ref{lem6}, we obtain
\begin{align*}
4\alpha
&=4\theta_{4}+3\theta_{3}+2\theta^{1}_{2}+2\theta^{2}_{2}+\theta_{1} \\[1pt]
&=2\theta^{1}_{2}+2\theta^{2}_{2}+\theta_{1}+o(1) \\[1pt]
&\leq 2+o(1).
\end{align*}
Thus, $\alpha\leq \frac{1}{2}+o(1)$. As $\alpha=\frac{e(\tilde{G})}{e(Q^{k}_{n})}$,
$$\mathrm{ex}(Q^{k}_{n},C_{6})=e(\tilde{G})=\alpha e(Q^{k}_{n})\leq \left(\frac{1}{2}+o(1)\right)e(Q^{k}_{n}).$$
\hfill $\blacksquare$
\vskip0.2cm

To sharpen the bound above, we introduce the auxiliary graph $H_{v}$.

\begin{defi}\label{defi3}
For each $v\in V(Q^{k}_{n})$, the auxiliary graph $H_{v}$ is the graph with vertex set $V(H_{v})=N_{Q^{k}_{n}}(v)\setminus N_{\tilde{G}}(v)$. Two vertices $u_{1},u_{2}\in V(H_{v})$ are adjacent in $H_{v}$ if and only if they have a common neighbor $w$ in $\tilde{G}$, where $w\neq v$.
\end{defi}

By Definition~\ref{defi3},
\begin{equation}
\sum_{v\in V(Q^{k}_{n})}v(H_{v})=\sum_{v\in V(Q^{k}_{n})}(2n-d_{\tilde{G}}(v)).
\end{equation}

Recall that $\theta^{1}_{2}$ denotes the proportion of 4-cycles $C$ in $Q^{k}_{n}$ such that $C\cap \tilde{G}$ contains exactly two adjacent edges of $\tilde{G}$. By Definition~\ref{defi3}, we have
\begin{equation}
\sum_{v\in V(Q^{k}_{n})}e(H_{v})=\theta^{1}_{2}n_{k}(C_{4}) \quad \text{for even } k\geq 4.
\end{equation}

For later use, we recall the classical Mantel's theorem.

\begin{lem}[\cite{Mantel-1907}]\label{lem7}
Every triangle-free graph on $n$ vertices has at most $\lfloor\frac{n^{2}}{4}\rfloor$ edges.
\end{lem}

The following theorem gives an improved upper bound on the generalized Tur\'{a}n number $\mathrm{ex}(Q^{k}_{n},C_{6})$.

\begin{theorem}\label{theorem3}
If $k$ is an even integer with $k\geq 4$, then $\mathrm{ex}(Q^{k}_{n},C_{6})\leq (\sqrt{2}-1+o(1))e(Q^{k}_{n}).$
\end{theorem}
\f {\bf Proof.}
For even $k\geq 4$, the number of 2-paths in $\tilde{G}$ that are contained in a 4-cycle in $Q^{k}_{n}$ is equal to
$(4\theta_{4}+2\theta_{3}+\theta^{1}_{2})n_{k}(C_{4})$, where $\theta_{4}, \theta_{3}$ and $\theta^{1}_{2}$ are defined in Lemma~\ref{lem5} and $n_{k}(C_{4})$ denotes the number of $4$-cycles in $Q^{k}_{n}$.

If $k\geq 6$, we have
\begin{align*}
(4\theta_4 + 2\theta_3 + \theta_2^1) n_k(C_4)
\ge \sum_{w \in V(\tilde{G})} \binom{d_{\tilde{G}}(w)}{2} - n v(Q_n^k).
\end{align*}
Substituting (14) into the inequality above and applying Lemma~\ref{lem6}, we obtain
\begin{equation}
\sum_{v\in V(Q^{k}_{n})}e(H_{v})\geq \sum_{w\in V(\tilde{G})}\binom{d_{\tilde{G}}(w)}{2}-nv(Q^{k}_{n})-o(n_{k}(C_{4})).
\end{equation}

Similarly, if $k=4$, we have
\begin{align*}
(4\theta_4 + 2\theta_3 + \theta_2^1) n_4(C_4)
&= \sum_{w \in V(\tilde{G})} \binom{d_{\tilde{G}}(w)}{2}.
\end{align*}
Thus,
$$\theta_2^1 n_4(C_4)\geq \sum_{w \in V(\tilde{G})} \binom{d_{\tilde{G}}(w)}{2} -nv(Q^{4}_{n})- (4\theta_4 + 2\theta_3) n_4(C_4).$$
Substituting (14) into the inequality above and applying Lemma~\ref{lem6}, we obtain
\begin{equation}
\sum_{v\in V(Q^{4}_{n})}e(H_{v})\geq \sum_{w\in V(\tilde{G})}\binom{d_{\tilde{G}}(w)}{2}-nv(Q^{4}_{n})-o(n_{4}(C_{4})).
\end{equation}

By Lemma~\ref{lem2}(2) and Lemma~\ref{lem3}(2), we have
\begin{equation}
n_{k}(C_{4})=\binom{n}{2}k^{n}<\frac{n^{2}}{2}k^{n} \quad \text{for even } k\geq 6,
\end{equation}
and
\begin{equation}
n_{4}(C_{4})=(2n^{2}-n)4^{n-1}
< \frac{n^{2}}{2}4^{n}.
\end{equation}

On the one hand, combining inequalities (15)-(18) yields
\begin{align}
\sum_{v\in V(Q^{k}_{n})} e(H_{v})
&\geq \sum_{w\in V(\tilde{G})} \binom{d_{\tilde{G}}(w)}{2} - n v(Q^{k}_{n}) - o\left( \frac{n^{2}}{2} k^{n} \right) \nonumber \\
&= \frac{1}{2} \sum_{w\in V(Q^{k}_{n})} d^{2}_{\tilde{G}}(w) - e(\tilde{G}) - n v(Q^{k}_{n}) - o\left( \frac{n^{2}}{2} k^{n} \right)
\end{align}
for all even $k\geq 4$.

On the other hand, since $\tilde{G}$ is 6-cycle-free, $H_{v}$ contains no triangles for any $v\in V(Q^{k}_{n})$. By Lemma~\ref{lem7}, we have
\begin{equation}
e(H_{v})\leq \frac{v^{2}(H_{v})}{4}.
\end{equation}
Combining (20) with (13) yields
\begin{align}
\sum_{v \in V(Q^{k}_{n})} e(H_v)
&\leq \sum_{w \in V(Q^{k}_{n})} \frac{\big(2n - d_{\tilde{G}}(w)\big)^2}{4} \nonumber \\
&= v(Q^{k}_{n}) n^2
   - 2n e(\tilde{G})
   + \frac{1}{4} \sum_{w \in V(Q^{k}_{n})} d_{\tilde{G}}^2(w). \label{eq:final}
\end{align}

Combining (19) with (21) and applying the Cauchy-Schwarz inequality, we have
\begin{align*}
v(Q^{k}_{n})n^{2}-2ne(\tilde{G})+\frac{1}{4}\sum_{w\in V(Q^{k}_{n})} d_{\tilde{G}}^2(w)
&\geq \frac{1}{2}\sum_{w\in V(\tilde{G})}d^{2}_{\tilde{G}}(w)-e(\tilde{G})-nv(Q^{k}_{n})-o\left( \frac{n^{2}}{2} k^{n} \right) \\
v(Q^{k}_{n})n^{2}-2ne(\tilde{G})
&\geq \frac{1}{4}\sum_{w\in V(\tilde{G})}d^{2}_{\tilde{G}}(w)-e(\tilde{G})-nv(Q^{k}_{n})-o\left( \frac{n^{2}}{2} k^{n} \right) \\
v(Q^{k}_{n})n^{2}-2ne(\tilde{G})
&\geq \frac{1}{4v(Q^{k}_{n})}\left(\sum_{w\in V(\tilde{G})}d_{\tilde{G}}(w)\right)^{2}-e(\tilde{G})-nv(Q^{k}_{n})-o\left( \frac{n^{2}}{2} k^{n} \right) .
\end{align*}
Rearranging terms gives
$$v(Q^{k}_{n})n^{2}-2ne(\tilde{G})-\frac{e^{2}(\tilde{G})}{v(Q^{k}_{n})}+e(\tilde{G})+nv(Q^{k}_{n})+o\left( \frac{n^{2}}{2} k^{n} \right)\geq 0.$$\\
Dividing both sides of the inequality by $v(Q^{k}_{n})n^{2}$ gives
$$1-\frac{2e(\tilde{G})}{v(Q^{k}_{n})n}-\frac{e^{2}(\tilde{G})}{v^{2}(Q^{k}_{n})n^{2}}+\frac{e(\tilde{G})}{v(Q^{k}_{n})n^{2}}+\frac{1}{n}+o\left(1\right)\geq 0.$$
Substituting $\alpha=\frac{e(\tilde{G})}{e(Q^{k}_{n})}=\frac{e(\tilde{G})}{nv(Q^{k}_{n})}$ into the above inequality, we obtain
$$1-2\alpha-\alpha^{2}+\frac{\alpha}{n}+\frac{1}{n}+o(1)\geq 0.$$

Since $0<\alpha<1$ and $n\rightarrow \infty$, we have $\frac{\alpha}{n}=o(1)$ and $\frac{1}{n}=o(1)$. Hence,
$$1-2\alpha-\alpha^{2}+o(1)\geq 0,$$
which implies $\alpha\leq \sqrt{2}-1+o(1)$. As $\alpha=\frac{e(\tilde{G})}{e(Q^{k}_{n})}$, we obtain
$$\mathrm{ex}(Q^{k}_{n},C_{6})=e(\tilde{G})=\alpha e(Q^{k}_{n})\leq (\sqrt{2}-1+o(1))e(Q^{k}_{n}).$$
\hfill $\blacksquare$
\vskip0.2cm

\subsection{4-Cycle-free subgraph of $Q^{k}_{n}$}
\label{sec:4-Cycle-free}
In this section, we focus on the case $\ell=2$. Let $\tilde{G}$ denote a 4-cycle-free spanning subgraph of $Q^{k}_{n}$ with maximum number of edges. The following theorem provides a general upper bound on the generalized Tur\'{a}n number $\mathrm{ex}(Q^{k}_{n},C_{4})$ that holds for all even $k\geq 4$.

\begin{theorem}\label{theorem4}
If $k$ is an even integer with $k\geq 4$, then $\mathrm{ex}(Q^{k}_{n},C_{4})\leq \frac{3}{4}e(Q^{k}_{n}).$
\end{theorem}
\f {\bf Proof.}
Since $\tilde{G}$ is $C_{4}$-free, we have $\theta_{4}=0$, where $\theta_{4}$ is defined in Lemma~\ref{lem5}. Let $\alpha=\frac{e(\tilde{G})}{e(Q^{k}_{n})}$. By Lemma~\ref{lem5}, we obtain
\begin{align*}
\alpha
&=\frac{1}{4}(3\theta_{3}+2\theta^{1}_{2}+2\theta^{2}_{2}+\theta_{1})\\
&\leq \frac{3}{4}(\theta_{4}+\theta_{3}+\theta^{1}_{2}+\theta^{2}_{2}+\theta_{1}+\theta_{0})\\
&=\frac{3}{4}.
\end{align*}
Consequently,
$$\mathrm{ex}(Q^{k}_{n},C_{4})=e(\tilde{G})=\alpha e(Q^{k}_{n})\leq \frac{3}{4}e(Q^{k}_{n}).$$
\hfill $\blacksquare$
\vskip0.2cm

By Corollary~\ref{cor1}, every 2-path in $Q^{4}_{n}$ is contained in a unique 4-cycle. This property is used to improve the general bound from Theorem~\ref{theorem4} for the special case $k=4$. The following lemma provides two fundamental facts about the structure of subcubes $Q_{3}$ in $Q^{4}_{n}$, which are useful in our subsequent arguments.

\begin{lem}\label{lem8}
Let $Q^{4}_{n}$ be the $4$-ary $n$-cube, and let $n_{4}(Q_{3})$ denote the number of subcubes $Q_{3}$ in $Q^{4}_{n}$. Then the following statements hold.
\begin{enumerate}
\item [{\rm (1)}] Every $K_{1,3}$ in $Q^{4}_{n}$ is contained in a unique subcube $Q_{3}$.
\item [{\rm (2)}] $n_{4}(Q_{3})=\frac{n(n-1)(2n-1)}{3}4^{n-1}$.
\end{enumerate}
\end{lem}
\f {\bf Proof.}
(1) Let $K_{1,3}$ be centered at $u$ with leaves $v,w,x$, and let $u=u_{1}u_{2}\cdots u_{n}$. Without loss of generality, assume that $u$ and $v$ differ in the $i$-th coordinate, $u$ and $w$ differ in the $j$-th coordinate, and $u$ and $x$ differ in the $s$-th coordinate, where $i,j,s\in [n]$. By Definition~\ref{defi1}, at most two of $i,j$ and $s$ are equal. We distinguish the following two cases.

{\bf Case 1.} Exactly two of $i,j,s$ are equal.

Without loss of generality, assume that $1\leq i=j < s\leq n$, and denote $v=u_{1}\cdots (u_{i}+1)\cdots u_{s}\cdots u_{n}$, $w=u_{1}\cdots (u_{i}-1)\cdots u_{s}\cdots u_{n}$ and $x=u_{1}\cdots u_{i}\cdots (u_{s}+\varepsilon)\cdots u_{n}$ for some $\varepsilon\in\{+1,-1\}$. By Corollary~\ref{cor1}, the 2-path $(v,u,w)$ is contained in a unique 4-cycle $(v,u,w,a,v)$ with $a=u_{1}\cdots (u_{i}+2)\cdots u_{s}\cdots u_{n}$; the 2-path $(v,u,x)$ is contained in a unique 4-cycle $(v,u,x,b,v)$ with $b=u_{1}\cdots (u_{i}+1)\cdots (u_{s}+\varepsilon)\cdots u_{n}$; the 2-path $(w,u,x)$ is contained in a unique 4-cycle $(w,u,x,c,w)$ with $c=u_{1}\cdots (u_{i}-1)\cdots (u_{s}+\varepsilon)\cdots u_{n}$. Then $(a,v,b)$, $(b,x,c)$, $(a,w,c)$ are three 2-paths in $Q^{4}_{n}$. By Corollary~\ref{cor1}, the 2-path $(a,v,b)$ is contained in a unique 4-cycle $(a,v,b,d,a)$; the 2-path $(b,x,c)$ is contained in a unique 4-cycle $(b,x,c,d,b)$; the 2-path $(a,w,c)$ is contained in a unique 4-cycle $(a,w,c,d,a)$, where $d=u_{1}\cdots (u_{i}+2)\cdots (u_{s}+\varepsilon)\cdots u_{n}$. Let $V=\{u,v,w,x,a,b,c,d\}$. Then the induced subgraph $Q^{4}_{n}[V]\cong Q_{3}$. Therefore, $K_{1,3}$ is contained in a unique subcube $Q_{3}$.

{\bf Case 2.} $i,j,s$ are pairwise distinct.

Without loss of generality, assume that $1\leq i < j < s\leq n$, and denote $v=u_{1}\cdots (u_{i}+\varepsilon_{1})\cdots u_{j}\cdots u_{s}\cdots u_{n}$, $w=u_{1}\cdots u_{i}\cdots (u_{j}+\varepsilon_{2})\cdots u_{s}\cdots u_{n}$ and $x=u_{1}\cdots u_{i}\cdots u_{j}\cdots (u_{s}+\varepsilon_{3})\cdots u_{n}$, where $\varepsilon_{1},\varepsilon_{2},\varepsilon_{3}\in\{+1,-1\}$. By Corollary~\ref{cor1}, the 2-path $(v,u,w)$ is contained in a unique 4-cycle $(v,u,w,a,v)$ with $a=u_{1}\cdots (u_{i}+\varepsilon_{1})\cdots (u_{j}+\varepsilon_{2})\cdots u_{s}\cdots u_{n}$; the 2-path $(v,u,x)$ is contained in a unique 4-cycle $(v,u,x,b,v)$ with $b=u_{1}\cdots (u_{i}+\varepsilon_{1})\cdots u_{j}\cdots (u_{s}+\varepsilon_{3})\cdots u_{n}$; the 2-path $(w,u,x)$ is contained in a unique 4-cycle $(w,u,x,c,w)$ with $c=u_{1}\cdots u_{i}\cdots (u_{j}+\varepsilon_{2})\cdots (u_{s}+\varepsilon_{3})\cdots u_{n}$. Then $(a,v,b)$, $(b,x,c)$, $(a,w,c)$ are three 2-paths in $Q^{4}_{n}$. By Corollary~\ref{cor1}, the 2-path $(a,v,b)$ is contained in a unique 4-cycle $(a,v,b,d,a)$; the 2-path $(b,x,c)$ is contained in a unique 4-cycle $(b,x,c,d,b)$; the 2-path $(a,w,c)$ is contained in a unique 4-cycle $(a,w,c,d,a)$, where $d=u_{1}\cdots (u_{i}+\varepsilon_{1})\cdots (u_{j}+\varepsilon_{2})\cdots (u_{s}+\varepsilon_{3})\cdots u_{n}$. Let $V=\{u,v,w,x,a,b,c,d\}$. Then the induced subgraph $Q^{4}_{n}[V]\cong Q_{3}$. Therefore, $K_{1,3}$ is contained in a unique subcube $Q_{3}$.

\medskip


(2) Since $d_{Q^{4}_{n}}(u)=2n$, the number of $K_{1,3}$ centered at $u$ is $\binom{2n}{3}$ for any vertex $u\in V(Q^{4}_{n})$. Clearly, every $Q_{3}$ contains exactly eight $K_{1,3}$. By part (1), every $K_{1,3}$ is contained in a unique subcube $Q_{3}$. Double-counting the pairs $(K_{1,3},Q_{3})$ such that $K_{1,3}$ is a subgraph of the $Q_{3}$ gives
$$8n_{4}(Q_{3})=\binom{2n}{3}v(Q^{4}_{n}).$$
Hence
$$n_{4}(Q_{3})
=\frac{n(n-1)(2n-1)}{3}4^{n-1}.$$
\hfill $\blacksquare$
\vskip0.2cm

Let $Q_{3}\cap \tilde{G}$ be the graph with vertex set $V(Q_{3})\cap V(\tilde{G})$ and edge set $E(Q_{3})\cap E(\tilde{G})$. The following lemma gives an upper bound on the number of vertices of degree 3 in $Q_{3}\cap \tilde{G}$.

\begin{lem}\label{lem9}
For every subcube $Q_{3}$ in $Q^{4}_{n}$, the graph $Q_{3}\cap \tilde{G}$ contains at most two vertices of degree $3$. Moreover, if $Q_{3}\cap \tilde{G}$ contains exactly two vertices of degree $3$, then they are either at distance $3$ or adjacent in $Q_{3}$.
\end{lem}
\f {\bf Proof.}
If $Q_{3}\cap \tilde{G}$ contains no vertex of degree 3, then the statement holds trivially. Otherwise, let $u$ be a vertex of $Q_{3}\cap \tilde{G}$ with $d_{Q_{3}\cap \tilde{G}}(u)=3$, and let $N_{Q_{3}\cap \tilde{G}}(u)=\{v,w,x\}$. We consider the degrees of three neighbors of $u$ in $Q_{3}\cap \tilde{G}$ in the following two cases.

{\bf Case 1.} No neighbor of $u$ in $Q_{3}\cap \tilde{G}$ has degree 3.

\begin{figure}[h!]
\begin{center}
\includegraphics[height=6cm,width=15cm]{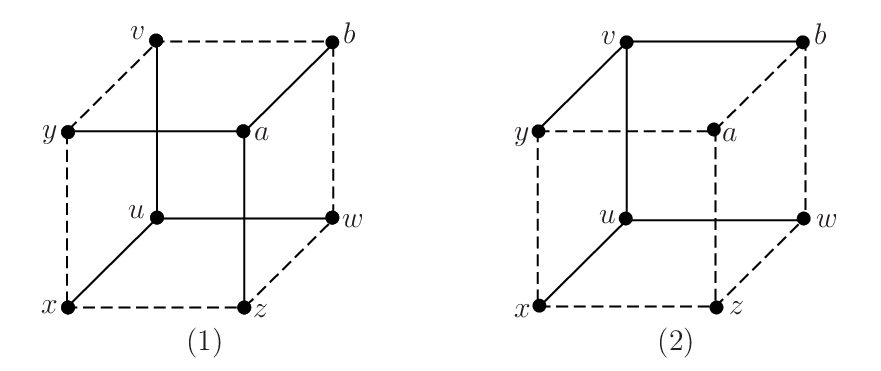}
\end{center}
\begin{center}
\vspace{-0.6cm}
\caption{ The two possible configurations of $Q_{3}\cap \tilde{G}$.}\label{F2}
\vspace{-0.6cm}
\end{center}
\end{figure}

In this case, no vertex at distance 2 from $u$ in $Q_{3}$ can have degree 3; otherwise $\tilde{G}$ contains a 4-cycle, contradicting the assumption that $\tilde{G}$ is 4-cycle-free. Thus, the only possible vertices of degree 3 in $Q_{3}\cap \tilde{G}$ are $u$ and the vertex at distance 3 from $u$. Hence there are at most two vertices of degree 3 in $Q_{3}\cap \tilde{G}$ (see Figure~\ref{F2}(1)).

{\bf Case 2.} Some neighbor of $u$ in $Q_{3}\cap \tilde{G}$ has degree 3.

Without loss of generality, assume that $v$ is a neighbor of $u$ with $d_{Q_{3}\cap \tilde{G}}(v)=3$. Then the other two neighbors $w$ and $x$ of $u$ cannot have degree 3; otherwise $\tilde{G}$ contains a 4-cycle, contradicting the assumption that $\tilde{G}$ is 4-cycle-free. Similarly, no vertex at distance 2 or 3 from $u$ can have degree 3. Therefore, the only vertices of degree 3 in $Q_{3}\cap \tilde{G}$ are $u$ and its neighbor $v$ (see Figure~\ref{F2}(2)).

Suppose $Q_{3}\cap \tilde{G}$ contains exactly two vertices of degree $3$. In Case 1, these two vertices are $u$ and the vertex at distance $3$ from $u$; in Case 2, they are $u$ and its neighbor $v$. Hence they are either at distance $3$ or adjacent in $Q_{3}$.
\hfill $\blacksquare$

In what follows, let $\sigma_{i}$ be the proportion of subcubes $Q_{3}$ in $Q^{4}_{n}$ such that $Q_{3}\cap \tilde{G}$ contains exactly $i$ vertices of degree $3$ for $i=0,1,2$. By Lemma~\ref{lem9}, the only possible values of $i$ are $0,1,2$. Hence
\begin{equation}
\sigma_{0}+\sigma_{1}+\sigma_{2}=1.
\end{equation}
\vskip0.2cm

\begin{lem}\label{lem10}
$2\sigma_{2}+\sigma_{1}=\dfrac{12}{n(n-1)(2n-1)4^{n}}\displaystyle\sum\limits_{w\in V(\tilde{G})}\binom{d_{\tilde{G}}(w)}{3}.$
\end{lem}
\f {\bf Proof.}
By Lemma~\ref{lem8}(1), every $K_{1,3}$ in $Q^{4}_{n}$ determines a unique subcube $Q_{3}$. Thus, by double counting the number of pairs $(w,Q_{3})$ with $w\in V(\tilde{G})$ and $Q_{3}\subseteq Q^{4}_{n}$ such that $d_{Q_{3}\cap \tilde{G}}(w)=3$, we obtain
$$(2\sigma_{2}+\sigma_{1})n_{4}(Q_{3})=\sum_{w\in V(\tilde{G})}\binom{d_{\tilde{G}}(w)}{3},$$\\
where $n_{4}(Q_{3})$ denotes the number of $Q_{3}$ in $Q^{4}_{n}$.
By Lemma~\ref{lem8}(2), we have
$$(2\sigma_{2}+\sigma_{1})\frac{n(n-1)(2n-1)}{3}4^{n-1}= \sum_{w\in V(\tilde{G})}\binom{d_{\tilde{G}}(w)}{3}.$$
Dividing both sides by $\frac{n(n-1)(2n-1)}{3}4^{n-1}$ yields
$$2\sigma_{2}+\sigma_{1}= \frac{12}{n(n-1)(2n-1)4^{n}}\sum_{w\in V(\tilde{G})}\binom{d_{\tilde{G}}(w)}{3}.$$
\hfill $\blacksquare$
\vskip0.2cm

The following theorem provides a preliminary upper bound on the generalized Tur\'{a}n number $\mathrm{ex}(Q^{4}_{n},C_{4})$.

\begin{theorem}\label{theorem5}
$\mathrm{ex}(Q^{4}_{n},C_{4})\leq (\dfrac{1}{4})^{\frac{1}{3}}(1+o(1))e(Q^{4}_{n}).$
\end{theorem}
\f {\bf Proof.}
As $0\leq \sigma_{i}\leq 1$ for $i=0,1$, by Lemma~\ref{lem10} and (22), we have
\begin{align}
\frac{12}{n(n-1)(2n-1)4^{n}}\sum_{w\in V(\tilde{G})}\binom{d_{\tilde{G}}(w)}{3}
&=2\sigma_{2}+\sigma_{1} \notag \\[-10pt]
&=2-2\sigma_{0}-\sigma_{1} \notag \\[4pt]
&\leq 2.  \label{eq:sigma_bound}
\end{align}
Let $\alpha=\frac{e(\tilde{G})}{e(Q^{4}_{n})}=\frac{e(\tilde{G})}{n4^{n}}$. Applying Jensen's inequality yields
\begin{align}
\sum_{w\in V(\tilde{G})}\binom{d_{\tilde{G}}(w)}{3}
&\geq 4^{n}\binom{\frac{1}{4^{n}}\sum_{w\in V(\tilde{G})}d_{\tilde{G}}(w)}{3} \notag\\
&=4^{n}\binom{2\alpha n}{3}.  \label{eq:jensen_result}
\end{align}
Combining (23) with (24), we obtain
$$\frac{12}{n(n-1)(2n-1)4^{n}}4^{n}\binom{2\alpha n}{3} \leq2.$$
Rearranging terms yields
$$8\alpha^{3}\leq 2+\frac{12\alpha^{2}-3}{n}+\frac{1-4\alpha}{n^{2}}.$$
Since $n\rightarrow\infty$ and $0<\alpha<1$, we have $\frac{12\alpha^{2}-3}{n}=o(1)$ and $\frac{1-4\alpha}{n^{2}}=o(1)$. Hence
$$8\alpha^{3}\leq 2+o(1),$$
which implies
$$\alpha\leq (\frac{1}{4})^{\frac{1}{3}}(1+o(1)).$$
As $\alpha=\frac{e(\tilde{G})}{e(Q^{4}_{n})}$, we have
$$\mathrm{ex}(Q^{4}_{n},C_{4})=e(\tilde{G})=\alpha e(Q^{4}_{n})\leq (\frac{1}{4})^{\frac{1}{3}}(1+o(1))e(Q^{4}_{n}).$$
\hfill $\blacksquare$
\vskip0.2cm

To improve the above bound, we introduce the following definitions:

For each $v\in V(Q^{4}_{n})$, let $H_{v}$ denote the auxiliary graph defined in Definition~\ref{defi3}. For $x,y\in V(H_{v})$, the pair $\{x,y\}$ is said to be \textit{blue} if $xy\in E(H_{v})$. Given a subcube $Q_{3}$ of $Q^{4}_{n}$, a vertex $v\in V(Q_{3})$ is called \textit{admissible} if its three neighbors $x,y,z$ in $Q_{3}$ all belong to $V(H_{v})$, and exactly $0$ or $3$ of the pairs $\{x,y\}$, $\{x,z\}$ and $\{y,z\}$ are blue.
A $3$-vertex subset $\{v_{1}, v_{2}, v_{3}\}\subseteq V(G)$ is called \textit{monochromatic} if either $v_{i}v_{j}\in E(G)$ for all distinct $i,j\in [3]$, or $v_{i}v_{j}\notin E(G)$ for all distinct $i,j\in [3]$.

The following classical result of Goodman~\cite{Goodman-1959} provides a lower bound on the number of monochromatic triangles in a graph.

\begin{lem}[\cite{Goodman-1959}]\label{lem11}
Let $G$ be a graph on $t$ vertices with edge density $p$. Then the number of monochromatic triangles in $G$ is at least $(1-3p+3p^{2})\binom{t}{3}$.
\end{lem}

By Lemma~\ref{lem9}, for every subcube $Q_{3}$ in $Q^{4}_{n}$, $Q_{3}\cap \tilde{G}$ has at most two vertices of degree $3$. We analyze the number of admissible vertices of $Q_{3}\cap \tilde{G}$ by considering the following three cases.

\begin{lem}\label{lem12}
Let $Q_{3}$ be a subcube of $Q^{4}_{n}$. Then the following hold for $Q_{3}\cap \tilde{G}$.
\begin{enumerate}
\item [{\rm (1)}] If $Q_{3}\cap \tilde{G}$ contains two vertices of degree $3$, then it contains no admissible vertices.
\item [{\rm (2)}] If $Q_{3}\cap \tilde{G}$ contains one vertex of degree $3$, then it contains at most one admissible vertex.
\item [{\rm (3)}] If $Q_{3}\cap \tilde{G}$ contains no vertex of degree $3$, then it contains at most eight admissible vertices.
\end{enumerate}
\end{lem}
\f {\bf Proof.}
Note that every admissible vertex has degree 0 in $Q_{3}\cap \tilde{G}$.

(1) By Lemma~\ref{lem9}, if $Q_{3}\cap \tilde{G}$ contains two vertices of degree $3$, then they are either at distance 3 or adjacent in $Q_{3}$. We consider these two cases separately.

{\bf Case 1.} The two vertices of degree 3 in $Q_{3}\cap \tilde{G}$ are at distance 3 in $Q_{3}$.

\begin{figure}[h!]
\begin{center}
\includegraphics[height=6cm,width=16cm]{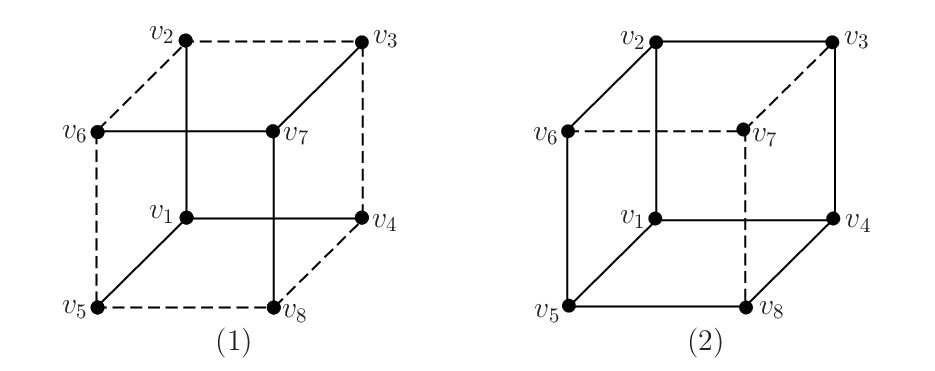}
\end{center}
\begin{center}
\vspace{-0.6cm}
\caption{ The illustration of Lemma~\ref{lem12}(1).}\label{F3}
\vspace{-0.6cm}
\end{center}
\end{figure}

Let $v_{1}$ and $v_{7}$ be these two vertices, with $d_{Q_{3}\cap \tilde{G}}(v_{1})=d_{Q_{3}\cap \tilde{G}}(v_{7})=3$ (see Figure~\ref{F3}(1)). Then each of the remaining six vertices in $Q_{3}$ has degree at least 1 in $Q_{3}\cap \tilde{G}$. Since admissible vertices have degree 0 in $Q_{3}\cap \tilde{G}$, $Q_{3}\cap \tilde{G}$ contains no admissible vertices.

{\bf Case 2.} The two vertices of degree 3 in $Q_{3}\cap \tilde{G}$ are adjacent in $Q_{3}$.

Let $v_{1}$ and $v_{2}$ be these two vertices, with $d_{Q_{3}\cap \tilde{G}}(v_{1})=d_{Q_{3}\cap \tilde{G}}(v_{2})=3$. Then at most two vertices, denoted by $v_{7}$ and $v_{8}$, have degree 0 in $Q_{3}\cap \tilde{G}$. By the symmetry of $Q_{3}$, assume that $v_{7}$ is admissible. Then the number of blue pairs among $\{v_{3},v_{6}\}$, $\{v_{6},v_{8}\}$ and $\{v_{3},v_{8}\}$ is either $0$ or $3$. Since the pair $\{v_{3},v_{6}\}$ is blue, all three pairs are blue, contradicting the assumption that $\tilde{G}$ is 4-cycle-free (see Figure~\ref{F3}(2)). Therefore, $Q_{3}\cap \tilde{G}$ contains no admissible vertices.

(2) Let $v_{1}$ be the unique vertex in $Q_{3}\cap \tilde{G}$ with $d_{Q_{3}\cap \tilde{G}}(v_{1})=3$. Then any vertex of degree 0 in $Q_{3}\cap \tilde{G}$ must be either at distance 2 or at distance 3 from $v_{1}$ in $Q_{3}$. If $Q_{3}\cap \tilde{G}$ contains no admissible vertex, the conclusion holds trivially. Otherwise, we consider these two cases separately.

{\bf Case 1.} Some vertex at distance $2$ from $v_{1}$ in $Q_{3}$ is admissible.

\begin{figure}[h!]
\begin{center}
\includegraphics[height=6cm,width=16cm]{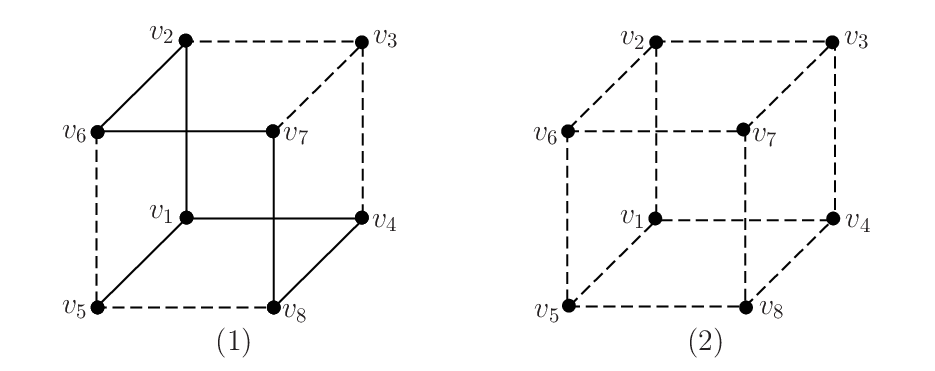}
\end{center}
\begin{center}
\vspace{-0.6cm}
\caption{ The illustration of Lemma~\ref{lem12}(2) and (3).}\label{F4}
\vspace{-0.6cm}
\end{center}
\end{figure}

Without loss of generality, let $v_{3}$ be such an admissible vertex. Then the number of blue pairs among $\{v_{2},v_{4}\}$, $\{v_{2},v_{7}\}$ and $\{v_{4},v_{7}\}$ is either $0$ or $3$. Since the pair $\{v_{2},v_{4}\}$ is blue, all three pairs are blue (see Figure~\ref{F4}(1)). Then no vertex other than $v_{3}$ has degree 0. Thus, $v_{3}$ is the unique admissible vertex.

{\bf Case 2.} The vertex at distance $3$ from $v_{1}$ in $Q_{3}$ is admissible.

Let $v_{7}$ be such an admissible vertex. Then the number of blue pairs among $\{v_{3},v_{6}\}$, $\{v_{6},v_{8}\}$ and $\{v_{3},v_{8}\}$ is either $0$ or $3$. If all three pairs are blue, this contradicts the assumption that $\tilde{G}$ is 4-cycle-free. Therefore, the number of blue pairs among them must be 0. It is possible that $v_{3}$, $v_{6}$ and $v_{8}$ each have degree 0. We claim that none of $v_{3},v_{6},v_{8}$ is admissible. Suppose to the contrary that one of them is admissible. By symmetry, we may assume that $v_{3}$ is admissible. Then the number of blue pairs among $\{v_{2},v_{4}\}$, $\{v_{2},v_{7}\}$ and $\{v_{4},v_{7}\}$ is either $0$ or $3$. Since the pair $\{v_{2},v_{4}\}$ is blue, all three pairs are blue, contradicting the assumption that $v_{7}$ is admissible (see Figure~\ref{F4}(1)). Therefore, $v_{7}$ is the unique admissible vertex.

(3) If $Q_{3}\cap \tilde{G}$ contains no vertex of degree 3, then it contains at most eight admissible vertices. This bound is trivial since $Q_3$ has exactly eight vertices (see Figure~\ref{F4}(2)).
\hfill $\blacksquare$
\vskip0.2cm

By Definition~\ref{defi3}, $v(H_{v})=2n-d_{\tilde{G}}(v)$. Let $\overline{d}_{\tilde{G}}(v):=v(H_{v})$, and let $p$ denote the edge density of $H_{v}$, i.e., $e(H_{v})=p\binom{\overline{d}_{\tilde{G}}(v)}{2}$. Recall that $\sigma_{0}$ and $\sigma_{1}$ denote the proportions of $Q_{3}$ in $Q^{4}_{n}$ such that $Q_{3}\cap \tilde{G}$ contains no vertex of degree $3$ and exactly one vertex of degree $3$, respectively.

\begin{lem}\label{lem13}
$\sigma_{1}+2\sigma_{0}\geq \dfrac{1}{16n_{4}(Q_{3})}\displaystyle\sum\limits_{v\in V(\tilde{G})}\binom{\overline{d}_{\tilde{G}}(v)}{3}.$
\end{lem}
\f {\bf Proof.}
On the one hand, by Lemma~\ref{lem12}, the total number of admissible vertices in all subcubes $Q_{3}$ of $Q^{4}_{n}$ is at most $(\sigma_{1}+8\sigma_{0})n_{4}(Q_{3})$, where $n_{4}(Q_{3})$ denotes the number of $Q_{3}$ in $Q^{4}_{n}$. On the other hand, for each $v\in V(\tilde{G})$, the graph $H_{v}$ has $\overline{d}_{\tilde{G}}(v)$ vertices and $p\binom{\overline{d}_{\tilde{G}}(v)}{2}$ edges. By Lemma~\ref{lem11}, $H_{v}$ contains at least $(1-3p+3p^{2})\binom{\overline{d}_{\tilde{G}}(v)}{3}$ monochromatic triangles. Since the number of admissible vertices in all subcubes $Q_{3}$ of $Q^{4}_{n}$ equals the total number of monochromatic triangles in $H_{v}$ over all $v\in V(\tilde{G})$, and $1-3p+3p^{2}\geq \frac{1}{4}$ for $0\leq p\leq 1$, we obtain
\begin{align}
(\sigma_1 + 8\sigma_0) n_4(Q_3)
&\ge (1 - 3p + 3p^2) \sum_{v \in V(\tilde{G})} \binom{\overline{d}_{\tilde{G}}(v)}{3} \notag \\
&\ge \frac{1}{4} \sum_{v \in V(\tilde{G})} \binom{\overline{d}_{\tilde{G}}(v)}{3}. \label{eq:admissible_bound}
\end{align}
Since $0\leq \sigma_{1}\leq 1$, we have
\begin{equation}
4\sigma_{1}+8\sigma_{0}\geq \sigma_{1}+8\sigma_{0}.
\end{equation}
Combining (25) with (26) yields
$$(4\sigma_{1}+8\sigma_{0})n_{4}(Q_{3})
\geq \frac{1}{4}\sum_{v\in V(\tilde{G})}\binom{\overline{d}_{\tilde{G}}(v)}{3}.$$
Dividing both sides by $4n_{4}(Q_{3})$, we obtain
$$\sigma_{1}+2\sigma_{0}
\geq \frac{1}{16n_{4}(Q_{3})}\sum_{v\in V(\tilde{G})}\binom{\overline{d}_{\tilde{G}}(v)}{3}.$$
\hfill $\blacksquare$
\vskip0.2cm

The following theorem improves the bound in Theorem~\ref{theorem5}.

\begin{theorem}\label{theorem6}
$\mathrm{ex}(Q^{4}_{n},C_{4})=\alpha e(Q^{4}_{n})\leq 0.628e(Q^{4}_{n}),$ where $\alpha$ satisfies $5\alpha^{3}+\alpha^{2}-\alpha-1\leq 0$.
\end{theorem}
\f {\bf Proof.}
Combining (23) with Lemma~\ref{lem13}, we have
\begin{align}
2
&= \frac{12}{n(n-1)(2n-1)4^{n}}\sum_{v\in V(\tilde{G})}\binom{d_{\tilde{G}}(v)}{3}
   + \sigma_{1}+2\sigma_{0} \nonumber \\
&\geq \frac{12}{n(n-1)(2n-1)4^{n}}\sum_{v\in V(\tilde{G})}\binom{d_{\tilde{G}}(v)}{3}
   + \frac{1}{16n_{4}(Q_{3})}\sum_{v\in V(\tilde{G})}\binom{\overline{d}_{\tilde{G}}(v)}{3} \nonumber \\
&= \frac{1}{n(n-1)(2n-1)4^{n+1}}
   \left(48\sum_{v\in V(\tilde{G})}\binom{d_{\tilde{G}}(v)}{3}
   + 3\sum_{v\in V(\tilde{G})}\binom{\overline{d}_{\tilde{G}}(v)}{3}\right). \label{eq:final}
\end{align}
Let $\alpha=\frac{e(\tilde{G})}{e(Q^{4}_{n})}=\frac{e(\tilde{G})}{n4^{n}}$. By Jensen's inequality and $\sum_{v\in V(\tilde{G})}\overline{d}_{\tilde{G}}(v) = 2(1-\alpha)n4^n$, we have
\begin{align}
\sum_{v\in V(\tilde{G})}\binom{\overline{d}_{\tilde{G}}(v)}{3}
&\geq 4^{n}\binom{\frac{1}{4^{n}}\sum_{v\in V(\tilde{G})}\overline{d}_{\tilde{G}}(v)}{3} \notag \\
&= 4^{n}\binom{2(1-\alpha)n}{3}. \label{eq:jensen_bound}
\end{align}
Substituting (24) and (28) into (27), we obtain
$$2\geq \frac{1}{n(n-1)(2n-1)4}\left(48\binom{2\alpha n}{3}+3\binom{2(1-\alpha) n}{3}\right).$$
Rearranging terms yields
$$(10\alpha^{3}+2\alpha^{2}-2\alpha-2)n^{3}+(-17\alpha^{2}+2\alpha+3)n^{2}+(5\alpha-1)n\leq 0.$$
Since $n\rightarrow\infty$, we have
$$10\alpha^{3}+2\alpha^{2}-2\alpha-2\leq 0,$$
or equivalently,
$$5\alpha^{3}+\alpha^{2}-\alpha-1\leq 0.$$
This implies $\alpha\leq 0.628$.
Since $\alpha=\frac{e(\tilde{G})}{e(Q^{4}_{n})}$, we obtain
$$\mathrm{ex}(Q^{4}_{n},C_{4})=\alpha e(Q^{4}_{n})\leq 0.628e(Q^{4}_{n}).$$
\hfill $\blacksquare$
\vskip0.2cm

\section{Concluding remarks}
\label{sec:remarks}

The generalized Tur\'{a}n number $\mathrm{ex}(G,H)$, a natural extension of the classical Tur\'{a}n number, is defined as the maximum number of edges in an $H$-free spanning subgraph of $G$. For each integer $\ell\geq 2$, let $C_{2\ell}$ denote the even cycle of length $2\ell$. In this paper, we study the generalized Tur\'{a}n number $\mathrm{ex}(Q^{k}_{n}, C_{2\ell})$ for all even integers $k\geq 4$ and all integers $\ell\geq 2$, thereby extending the known result for the hypercube $Q_{n}$ [J. Graph Theory 16 (1992) 273--286] from $k=2$ to $k\geq 4$. In future work, we plan to investigate the extremal problem $\mathrm{ex}(Q^{k}_{n}, C_{2\ell})$ for odd $\ell\geq 5$, which appears to be a challenging problem.


\end{document}